\documentclass[draft,12pt]{article}

\usepackage[latin1]{inputenc}
\usepackage{amsmath,amssymb}
\usepackage{latexsym}
\usepackage{amsthm}

\newtheorem{theorem}{Theorem}[section]
\newtheorem{corollary}[theorem]{Corollary}
\newtheorem{lemma}[theorem]{Lemma}
\newtheorem{proposition}[theorem]{Proposition}
\newtheorem{definition}[theorem]{Definition}

\newtheorem{remark}[theorem]{Remark}
\newtheorem{question}[theorem]{Question}

\numberwithin{equation}{section}

\def\qed{{\hfill $\square$ \bigskip}}

\def\square{{\vcenter{\vbox{\hrule height.3pt
        \hbox{\vrule width.3pt height5pt \kern5pt
           \vrule width.3pt}
        \hrule height.3pt}}}}

\def\sS {{\cal S}}

\def\bone{{\bf 1}}

\def\wt{\widetilde}

\def\P{{\mathbb P}}

\def\lam{{\lambda}}

\def\bee{\begin{equation}}
\def\bet{\begin{theorem}}
\def\bep{\begin{proposition}}
\def\bef{\begin{proof}}
\def\bel{\begin{lemma}}
\def\bec{\begin{corollary}}
\def\bed{\begin{definition}}
\def\ber{\begin{remark}}
\def\bequ{\begin{question}}
\def\eee{\end{equation}}
\def\eet{\end{theorem}}
\def\eep{\end{proposition}}
\def\eef{\end{proof}}
\def\eel{\end{lemma}}
\def\eec{\end{corollary}}
\def\eed{\end{definition}}
\def\eer{\end{remark}}
\def\eequ{\end{question}}

\def\R{{\mathbb R}}

\def\P{{\mathbb P}}

\def\lam{{\lambda}}

\def\al{{\alpha}}

\def\qed{{\hfill $\square$ \bigskip}}

\def\eps{\varepsilon}

\def\wt{\widetilde}

\def\ni{\noindent }
\def\ms{\medskip}

\def \half {{{1/ 2}}}

\def\square{{\vcenter{\vbox{\hrule height.3pt
        \hbox{\vrule width.3pt height5pt \kern5pt
           \vrule width.3pt}
        \hrule height.3pt}}}}

\def\tlint{{- \kern-0.85em \int \kern-0.2em}}  % for textstyle
\def\dlint{{- \kern-1.05em \int \kern-0.4em}}  % for displays

\def\sS {{\cal S}}

\begin{document}

\title{The supremum of Brownian local times on H\"older curves, II\\
{$\empty$}\\
Le supremum du temps locaux d'un mouvement Brownien sur les courbes Holderiennes, II}
\author{Richard F. Bass and Krzysztof Burdzy
\footnote{\rm Research of KB partially supported by Simons Foundation grant 928958.}
}

%\date{\today}

\maketitle

\begin{abstract}  
\noindent {\it Abstract:} 
For $f: [0,1]\to \R$, we consider $L^f_t$,
the local time of space-time Brownian motion on the curve
$f$.
Let $\sS_\al$ be the class of all functions whose H\"older norm of order $\al$
is less than or equal to 1. We show that the supremum of $L^f_1$ over
$f$ in $\sS_\al$ is finite if $\al>\frac12$.

\ms
\ni {\it  Abstrait:}
Soit $W_t$ un mouvement brownien et soit $L^f_t$ le temps local
du processus $(t,W_t)$ pour le  courbe $f:[0,1]\to \R$. 
Soit $\sS_\al$ la classe de toutes fonctions telle que la norme
holderienne du ordre $\al$ est moins de 1. Nous d\'emontrons que
$\sup_{f\in {\sS_\al}} L^f_1 <\infty$ p.s. si $\al>\frac12$.

\end{abstract}

\section{Introduction}

 The main  claim of \cite{BaBu1}  was that the supremum of Brownian local 
times over all $\alpha$-H\"older curves is finite if $\alpha >1/2$ 
(see Theorem \ref{T3.6} below for the precise statement). An error in the 
proof was pointed out to us by A.\ Vatamanelu; however we were able to 
establish the claim for $\al>5/6$ in \cite{BaBu2}.  
The purpose of this note is to prove 
the original claim from \cite{BaBu1}, that finiteness of the supremum indeed 
holds for all $\alpha\in(1/2,1]$. 
We also showed in \cite{BaBu1} that 
$\al=1/2$ is the critical value;
see Theorem 3.8 of that paper for the precise statement.

Let $W_t$ be one-dimensional Brownian motion and let $f:[0,1]\to \R$
be a H\"older continuous function. There are a number of equivalent ways to define $L^f_t$, the 
local time  of $W_t$ along the curve $f$, one of which is
 as the limit in probability of
$$\frac{1}{2\eps}\int_0^t \bone_{(f(s)-\eps, f(s)+\eps)}(W_s)\, ds$$
as $\eps\to 0$. See \cite[Sect. 2]{BaBu1} for a discussion of other ways of defining
$L^f_t$.
Let
$$\sS_\al=\{f: \sup_{0\leq t\leq 1}|f(t)|\leq 1,  
|f(s)-f(t)|\leq |s-t|^\al\mbox{ if }s,t\leq 1\}.$$

Our main result in this paper is the following.

\bet\label{T3.6}  For any $\al \in (1/2,1]$,
there exists $\wt L^f_t$ such that
\item{(i)} for each $f\in \sS_\al$, we have
$\wt L^f_t=L^f_t$ for all $t$, a.s.,   
\item{(ii)} with probability one, $ f\to \wt L^f_1$ is a continuous
map on $\sS_\al$ with respect to the supremum norm, and
\item{(iii)} with probability one, $\sup_{f\in \sS_\al} \wt L^f_1 <\infty$.
\eet

The interest in Theorem \ref{T3.6} has several sources. One is that the metric
entropy of $\sS_\al$ is far too large for chaining arguments to work; nevertheless the supremum is finite a.s. 
Another is the work of Holden and Sheffield \cite{HS} on scaling limits of the Schelling model, where they used some of the techniques in \cite{BaBu1} to analyze 
local times of random fields over Lipschitz surfaces.

In the interests of space we present only the changes needed to \cite{BaBu1}
to prove our result and refer to the original paper for the unchanged part of
the proof.

\section{The finiteness of the supremum}

Let $W_t$ be a Brownian motion.
A key ingredient in our proof is Lemma 3.1 of \cite{BaBu1}.
The proof there is correct; the error in \cite{BaBu1} was in how this lemma was applied further on.

We replace Propositions 3.2 and 3.3 in \cite{BaBu1} by the following.

Consider an integer $N>0$. 
For $0\leq \ell\leq N,
-N^\al-1\leq m\leq N^\al$, 
let $R_{\ell m}$ be the rectangle defined by
$$R_{\ell m}= [\ell/N, (\ell+1)/N]\times [m/N^\al, (m+1)/N^\al].$$

\bep\label{P3.2}
Let $\al \in (1/2,1]$ and $\eps\in (0,1/16)$.
There exist $c_1, c_2$, and $c_3$ such that:
\item{(i)} there exists a set $D_N$ with $\P(D_N)\leq c_1\exp(-c_2N^{\eps/2})$;
\item{(ii)} if $\omega\notin D_N$ and $f\in \sS_\al$, 
then there are at most $c_3 N^{(1/2)+2\eps}$ 
rectangles $R_{\ell m}$ in $[0,1]\times [-1,1]$ which contain 
both a point of the graph of $f$ and a point of the graph of $W_t(\omega)$.
\eep

\proof
Let $M=\lfloor N^\eps\rfloor$ and set
$$Q_{ik} = [i/M, (i+1)/M]\times [k/M^\half,(k+1)/M^\half],$$
for $0\leq i\leq M$ and $ -M^\half-1 \leq k\leq M^\half$.
Let $J=\lceil N/M \rceil$.

Let 
\begin{align*}I_{ikj}=\{ \exists t\in [i/M+(j-1)/N,&\, i/M+j/N]:\\ 
&k/M^\half\leq W_t\leq (k+1)/M^\half \},
\end{align*}
$$A_{ik}=\sum_{j=1}^J \bone_{I_{ikj}},$$
and
$$C_{ik}=\{A_{ik}\geq J^{(1/2)+\eps}\}.$$
By Lemma 3.1 of \cite{BaBu1}  with  $\lam = J^\eps$ and the Markov property
applied at $i/M$
we have $\P(C_{ik})\leq c_4\exp(-c_5 J^{\eps})$.

There are at most $c_6M^{3/2}$ rectangles $Q_{ik}$, so if
$D_N=\cup_{i,k} C_{ik}$, where 
$0\leq i\leq M$ and $ -M^\half-1 \leq k\leq M^\half$,
then
$$\P(D_N)\leq c_7 M^{3/2} \exp(-c_5 J^{\eps})
\leq c_7  \exp(-c_8 N^{\eps/2}).$$

 Let $f$ be
any function in $\sS_\al$. If $f$ intersects $Q_{ik}$ for some $i$ and
$k$, then $f$ might intersect $Q_{i,k-1}$ and $Q_{i,k+1}$.
But because $f\in\sS_\al$ and $\al>1/2$, it cannot
intersect $Q_{ir}$ for any $r$ such that $|r-k|>1$. Therefore
$f$ can intersect at most $3(M+1)$ of the $Q_{ik}$. 

Now suppose $\omega\notin  D_N$.
Look at any one of the $Q_{ik}$ that $f$ intersects. Since $\omega\notin
D_N$, then there are at most $J^{(1/2)+\eps}$ integers $j$ that
are less than $J$ and for which the path of $W_t(\omega)$ intersects
$([i/M+(j-1)/N,i/M+j/N]\times[-1,1])\cap Q_{ik}$.
If $f$ intersects a rectangle $R_{\ell m}$,
then it can intersect a rectangle $R_{\ell r}$
only if $|r-m|\leq 1$, since $f\in\sS_\al$.
Therefore
 there are at most $3J^{(1/2)+\eps}$ rectangles
$R_{\ell m}$ contained
in $Q_{ik}$ which contain both a point of the graph of $f$
and a point of the graph of $W_t(\omega)$.

Since there are at most $3(M+1)$ rectangles $Q_{ik}$ which contain
a point of the graph of $f$, there are therefore at most
$$3(M+1) 3J^{(1/2)+\eps}\le c_9 N^{(1/2)+2\eps}$$
rectangles $R_{\ell m}$ that
contain both a point of the graph of $f$ and a point of the graph of $W_t(\omega)$.
\qed

Our present Proposition \ref{P3.2} is almost identical to Proposition 3.3 in \cite{BaBu1}, so the latter can be omitted.
 With this change, the remainder of \cite{BaBu1}, beyond Proposition 3.3, can proceed as in the original.

\medskip

\ni {\bf Richard F. Bass}\\
Department of Mathematics\\
University of Connecticut \\
Storrs, CT 06269-3009, USA\\
{\tt r.bass@uconn.edu} 
\ms

\ni {\bf Krzysztof Burdzy}\\
Department of Mathematics\\
University of Washington\\
Box 354350\\
Seattle, WA 98195-4350\\
{\tt  burdzy@uw.edu}


\begin{thebibliography}{99}



\bibitem{BaBu1}
R.F. Bass and K. Burdzy (2001). 
The supremum of Brownian local times on H\"older curves. 
{\sl Ann. I.H. Poincar\'e \bf  37}, 38 (2002) 627--642.

\bibitem{BaBu2}
R.F. Bass and K. Burdzy (2002). 
Erratum to ``The supremum of Brownian local times on H\"older curves.'' 
{\sl Ann. I.H. Poincar\'e \bf  38}, 799--800.


\bibitem{HS}
N. Holden and S. Sheffield (2020)
{\sl Probab. Th. rel. Fields \bf 176}, 219--292.


\end{thebibliography}
\end{document}